\documentclass[12pt]{amsart}

\usepackage{amsthm}
\theoremstyle{plain}
\newtheorem{theorem}{Theorem}[section]
\newtheorem{lemma}[theorem]{Lemma}
\newtheorem{proposition}[theorem]{Proposition}
\newtheorem{corollary}[theorem]{Corollary}

\theoremstyle{definition}
\newtheorem{definition}[theorem]{Definition}
\newtheorem{example}[theorem]{Example}

\theoremstyle{remark}
\newtheorem{remark}{Remark}[theorem]
\newtheorem*{claim}{Claim}

\usepackage{amssymb}
\usepackage{amsxtra}
\usepackage[mathscr]{eucal}
\usepackage{enumerate}
\usepackage[all]{xy}

\usepackage{times,mathptmx}

\def\lto{\longrightarrow}

\DeclareMathOperator{\Pic}{Pic}
\DeclareMathOperator{\Hilb}{Hilb}
\DeclareMathOperator{\Sym}{Sym}
\DeclareMathOperator{\trace}{trace}
\DeclareMathOperator{\td}{td}
\DeclareMathOperator{\id}{id}
\DeclareMathOperator{\Aut}{Aut}

\begin{document}
\title{A characterization of certain irreducible symplectic 4-folds}
\author{Yasunari Nagai}                     

\address{Yasunari Nagai, 
Graduate School of Mathematical Sciences, University of Tokyo, 
3-8-1 Komaba, Meguro, Tokyo 153-8914, Japan}
\email{nagai@ms.u-tokyo.ac.jp}

%\institute{Graduate School of Mathematical Sciences, 
%University of Tokyo, 3-8-1 Komaba, Meguro, Tokyo 153-8914, Japan
%\email{nagai@ms.u-tokyo.ac.jp}}


%
\thanks{\textit{Mathematics Subject Classification} (2000): Primary 14J32, Secondary 32Q20}
\date{Revised 21/09/2002}

\maketitle

\begin{abstract}
We give a characterization of irreducible symplectic fourfolds
which are given as Hilbert scheme of points on a K3 surface.
\end{abstract}
\section{Introduction}
\label{intro}

In the theory of the moduli problem of K3 surfaces, 
Kummer surfaces played a very important role. It is easy to 
characterize Kummer surfaces. 

\begin{proposition}\label{kum} \textup{(See \cite{B-P-V}, Chap. VIII \S 6)} 
Let $S$ be a K3 surface. If $S$ contains 16 disjoint $\mathbb P^1$'s 
$C_1,\dots ,C_{16}$ and $D=\sum C_i$ is $2$-divisible
in $\Pic (S)$ then $S$ is isomorphic to a Kummer surface.
\end{proposition}

The density of Kummer surfaces in the moduli space and
this characterization enable us to derive the Global Torelli Theorem
for arbitrary K3 surfaces from that for Kummer surfaces.

A higher dimensional analogue of a K3 surface is an irreducible 
symplectic manifold.

\begin{definition}
A compact K\"ahler manifold $X$ of dimension $2n$ 
is said to be \emph{irreducible symplectic} if 
the following conditions are satisfied.
\begin{enumerate}[{\rm (i)}]
\item $X$ admits a symplectic form, \emph{i.e.} there exists a $d$-closed 
holomorphic 2-form $\sigma$ such that $\sigma ^{\wedge n}$ is nowhere 
vanishing.
\item $h^0(X,\Omega ^2 _X)=1$, \emph{i.e.} any non-zero holomorphic 2-form
is the symplectic form up to constant.
\item $X$ is simply connected.
\end{enumerate}
An irreducible symplectic manifold is also called 
\emph{hyper-K\"ahler} in the literature (see \cite{Be,H}).
\end{definition}

It seems that the moduli behaviour of irreducible symplectic manifolds 
is similar to that of K3 surfaces. Although Namikawa recently found 
a counterexample to the Global Torelli Problem in higher 
dimensions \cite{Na}, one still believes that some kind of Global 
Torelli Theorem should hold, but even a convincing conjectural version of 
it is missing for the time being.

With a view towards the Global Torelli Problem for irreducible symplectic
manifolds, it is important to ask for some ``typical'' objects 
in the moduli spaces in question and to give their characterization. 
This question in general seems to be quite hard. It is natural 
to restrict ourselves to a special case as our first step.

The Hilbert scheme of points on a K3 surface is an example of 
irreducible symplectic manifold which is important and seems to be 
rather easy to handle, for it has a very explicit description, 
in particular in dimension four.

\begin{example}\label{hilb} \textup{(cf. \cite{Be})} 
Let $S$ be a smooth surface, $\Hilb ^n (S)$ the 
Hilbert scheme of 0-dimensional sub-schemes of length $n$ and 
$\Sym ^n (S)=S^n/\mathfrak S_n$ the $n$-th symmetric product of S. 
Beauville \cite{Be} showed that 
the natural morphism (\emph{Hilbert-Chow morphism})
\[
F : \Hilb ^n (S)\to \Sym ^n (S)
\]
is a crepant birational morphism and that if $S$ is a K3 surface, 
the Hilbert scheme $\Hilb ^n(S)$ is an irreducible 
symplectic manifold of dimension $2n$. In the case
$n=2$, the description of $F$ is quite easy.
The singular locus $\Sigma$ of $\Sym ^2(S)$ is isomorphic to $S$
and $\Sym ^2(S)$ is locally 
of the form $\mathbb C^2\times (A_1\mbox{ surface singularity})$ 
along $\Sigma$. It is easy 
to show that $F$ is simply the blowing-up of $\Sym ^2(S)$ along
$\Sigma$. Considering the action of $\mathfrak S_2$, we have the
following diagram
\begin{equation}\label{hilbdiagram}
\vcenter{\xymatrix{
Bl_{\Delta}(S\times S) \ar[r]^{\quad \widetilde F} \ar[d] 
  & S\times S \ar[d]\\
\Hilb ^2(S) \ar[r]^F & \Sym ^2(S)
}}
\end{equation}
where $\Delta $ is the diagonal of $S\times S$.
\end{example}

We give the following result as an analogy of
Proposition \ref{kum}.

\begin{theorem}\label{main}
Let $X$ be a projective irreducible symplectic fourfold. Assume that 
there exists a birational morphism $f:X\to Y$ which contracts 
an irreducible divisor $E$ to a surface $S\subset Y$ such that 
\begin{enumerate}[{\rm (i)}]
\item The general fibre of $f_{|E}:E\to S$ is isomorphic to $\mathbb P^1$,
\item $E$ is $2$-divisible in $\Pic (X)$, 
\item $E^4=192$.
\end{enumerate}
Then, $S$ is a K3 surface and 
$X$ is isomorphic to the Hilbert scheme $\Hilb ^2(S)$
of $S$.
\end{theorem}

\begin{remark}
If $X$ is deformation equivalent 
to some $\Hilb ^2(T)$ for a K3 surface $T$, the condition (iii) can be 
replaced by 
\begin{list}{(iii')}{\setlength{\leftmargin}{7pt}}
\item $q_X(E)=-8$.
\end{list}
where $q_X$ is the Beauville-Bogomolov form on $H^2(X,\mathbb{Z})$
(see \cite{Be,H}, see also Remark \ref{conv}).
\end{remark}

\begin{remark}
It would be natural to pose the following question: \emph{
If $E$ is an irreducible divisor on $X$ with $q_X(E)<0$, then 
there exist 
an irreducible symplectic fourfold $X'$ birational to $X$ and
a birational morphism $f:X'\to Y'$ which contracts the strict
transform of $E$ on $X'$ ?}  Clearly, the answer will be affirmative 
if every flop of symplectic 4-fold is a Mukai flop as conjectured, for 
the termination of flops for terminal fourfolds is already known.
\end{remark}

The next natural problem to consider is the density of the 
\emph{birational (bimeromorphic) models} 
of Hilbert schemes made from Kummer surfaces in the connected 
component of the moduli space containing an irreducible symplectic 
fourfold which is birational to  
the Hilbert scheme of some K3 surface. 
But even this seems to be rather hard question.

The rough idea to prove the theorem is to trace backward Beauville's 
proof of Example \ref{hilb}. It uses more or less elementary and standard 
techniques. It contains several ingredients. 
One is a numerical computation 
using Holomorphic Lefschetz theorem of Atiyah-Singer. Another is 
the decomposition theorem of K\"ahler manifolds with trivial 
first Chern class. 
The result of Wierzba \cite{W} on divisorial contractions of 
symplectic manifolds is also used in an essential way.
The remaining part consists of geometric arguments based 
on the geometry of K3 surfaces.

\paragraph{Notation.}
Through this paper we work with the following notation. 
Let $X$, $Y$, $f$, $E$ and $S$ be as in the theorem above. 
Theorems 1.4(ii) and 1.5 in \cite{W} imply that $E$ is a $\mathbb P^1$ 
bundle and $S$, 
which is the singular locus of $Y$, is a smooth surface with $K_S\sim 0$. 
Furthermore they infer that $Y$ is analytically locally 
isomorphic to $\mathbb C^2\times (\mbox{$A_1$ surface singularity})$ 
at each point of $S$. 
Put $D=\frac 12 E$ and take a double covering 
$p:\widetilde{X}\to X$ defined by $\mathscr{O}(D)$. 
Then $p$ is ramified at $\widetilde E\subset \widetilde X$
and $p(\widetilde E)=E$. 
$\widetilde X$ is smooth since $E$ is smooth and
we have the following diagram
\begin{equation}\label{modeldiagram}
\vcenter{\xymatrix{
\widetilde X \ar[r]^{\tilde f} \ar[d]_p & \widetilde Y \ar[d]^q\\
X \ar[r]^f & Y
}}
\end{equation}
where $\tilde f$ and $q$ is the Stein factorization of $f\circ p$.
In fact $\widetilde Y$ is a smooth fourfold, 
$\tilde f$ is the blowing-up of $\widetilde Y$ along a smooth
centre $\widetilde S$ with the exceptional divisor $\widetilde E$, and
$q_{|\widetilde S}:\widetilde S\overset{\sim}{\to} S$ is an
isomorphism. Note that $K_{\widetilde Y}\sim 0$, for $f$ is 
crepant and $q$ has no ramification divisor.

\begin{remark}
The projectivity assumption in Theorem \ref{main} is made to 
apply Wierzba's result in our argument. If Wierzba's description on 
divisorial contraction of symplectic manifolds is valid for 
non-projective ones, the projectivity assumption would not be 
necessary.
\end{remark}

\section{Geometry of $\widetilde Y$}
\label{sec:2}

In this section, we prove the following proposition.

\begin{proposition}\label{P1}
Under the assumption of Theorem \ref{main} and the notation above, 
$\widetilde Y$ is isomorphic to a product of two K3 surfaces.
\end{proposition}

Our strategy is to apply the following famous 
decomposition theorem to $\widetilde Y$.

\begin{theorem}[cf. \cite{Be}]\label{Bogomolov}
Let $Z$ be a compact K\"ahler manifold with $K_Z\sim 0$. Then 
there exists a finite \'etale covering $\widetilde Z\to Z$ such
that $\widetilde Z$ is isomorphic to a product of varieties of following
types
\begin{enumerate}[{\rm (i)}]
\item complex torus,
\item Calabi-Yau manifold \emph{i.e.} compact K\"ahler manifold $W$ for which 
$K_W\sim 0$, $h^i(W,\mathscr O_W)=0$ for $0<i<\dim W$, and $\pi _1(W)
=\{e\}$,
\item irreducible symplectic manifold.
\end{enumerate}
\end{theorem}

Thanks to this powerful theorem, Proposition \ref{P1} is reduced to 
the following

\begin{proposition}\label{quant}
{\rm (i)} $\pi _1(\widetilde Y)=\{e\}$.\quad
{\rm (ii)} $h^0(\widetilde Y, \Omega ^2 _{\widetilde Y})=2$.
\end{proposition}

\begin{proof}[Proof of\; Proposition \ref{quant}
$\Rightarrow$ Theorem \ref{P1}]
Applying Theorem \ref{Bogomolov} under the condition (i) 
of Proposition \ref{quant}, $\widetilde Y$ itself decomposes into 
a product of Calabi-Yau manifolds and  
irreducible symplectic manifolds. Since $\widetilde Y$ is of dimension $4$, 
a product of two K3 surfaces, a Calabi-Yau fourfold or a compact 
irreducible symplectic fourfold is possible. But (ii) of 
Proposition \ref{quant} asserts that the last two cases do 
not happen. 
\end{proof}

To compute these quantities from the condition (iii) of Theorem \ref{main}, 
we use Holomorphic Lefschetz formula of Atiyah-Singer \cite{A-S}.

\begin{definition}
Let $M$ be a compact complex manifold, and $g$ an automorphism of finite
order of $M$. The \emph{holomorphic Lefschetz number} $L_{hol}(g)$ is defined
by 
\[
L_{hol}(g)=\sum (-1)^p \trace (g^*:H^p(M,\mathscr O_M)) .
\]
\end{definition}

\begin{theorem}[Holomorphic Lefschetz formula, \cite{A-S}]\label{hollef}
As in the notation of the definition above. Assume further
that the fixed point set $M^g=\{x\in M\mid g(x)=x\}$ is smooth. 
Then the formula
\[
L_{hol}(g)
=\int _{M^g}
	\frac{\prod _{\theta} \mathcal U^{\theta} (N_{M^g/M}(\theta))
		\cdot \td (M^g)}
	{\det (1-(g_{N_{M^g/M}})^*)}
\]
holds, where $\td(M^g)$ denotes the Todd class of $M^g$, 
$N_{M^g/M}$ the normal bundle, 
$N_{M^g/M}(\theta)\subset N_{M^g/M}$ the eigen-sub-bundle of 
$(g_{N_{M^g/M}})^*$ with the eigenvalue $e^{i\theta}$, and
\[
\mathcal U^{\theta} (x_1,x_2,\dots)
=\left\{
	\prod _j \left(
		\frac{1-e^{-x_j-i\theta}}{1-e^{-i\theta}}
	\right)
\right\}^{-1} .
\]
\end{theorem}

The general formula itself is very complicated, but in our case 
the formula becomes easy to handle.

\begin{lemma}\label{appdhollef}
Notation as in {\rm \S 1}. 
Let $g$ be an involution of $\widetilde Y$ induced 
by $q:\widetilde Y\to Y$. Then $L_{hol}(g)\in \mathbb Z$
and we have
\[
L_{hol}(g)=\frac{1}{48}(c_2(T_{\widetilde S})
	+3 c_2(N_{\widetilde S/\widetilde Y})).
\]
In particular if $\widetilde S$ is a K3 surface,
\[
L_{hol}(g)=\frac 12
	+\frac 1{16} c_2(N_{\widetilde S/\widetilde Y})
\]
and if $S$ is an abelian surface, 
\[
L_{hol}(g)=\frac 1{16} c_2(N_{\widetilde S/\widetilde Y}).
\]
\end{lemma}

\begin{proof}
Since $g$ is an involution, eigenvalues of $g^*$
on each cohomology group must be $\pm 1$, therefore 
$\trace g^*\in \mathbb Z$, in particular $L_{hol}(g)\in \mathbb Z$. 

We apply Theorem \ref{hollef} under
$M=\widetilde Y$ and $M^g=\widetilde S$.
Since $g$ produces the two dimensional locus of $A_1$ singularities $S$, 
we have 
\[
(g_{N_{\widetilde S/\widetilde Y}})^* =
\begin{pmatrix}
-1 & 0\\
0 & -1
\end{pmatrix}.
\]
Therefore, $\det (1-(g_{N_{\widetilde S/\widetilde Y}})^*)=4$
and the only possible $\theta$ is $\theta =\pi$, \emph{i.e.} 
$N_{\widetilde S/\widetilde Y}(\pi)=N_{\widetilde S/\widetilde Y}$.
Since the rank of this bundle is 2, we have only to consider
$\mathcal U^{\pi}$ in $2$ variables. By definition
\[
\begin{aligned}
\mathcal U^{\pi}(x_1,x_2)
&= \frac 2{1+e^{-x_1}} \cdot \frac 2{1+e^{-x_2}}\\
&=1-\frac 12(x_1+x_2)+\frac 14 x_1x_2 +O(x_1,x_2)^3 .
\end{aligned}
\]
This implies 
\[
\mathcal U^{\pi}(N_{\widetilde S/\widetilde Y}(\pi))
=1-\frac 12 c_1(N_{\widetilde S/\widetilde Y})
+\frac 14 c_2(N_{\widetilde S/\widetilde Y}) .
\]
Note that $\td (\widetilde S)=1+
\frac 1{12} c_2 (T_{\widetilde S})$ because $K_{\widetilde S}
\sim 0$. Combining these, we get the formula. For the last assertion,
note that $\widetilde S$ is either a K3 surface or a complex $2$-torus by 
the Enriques-Kodaira's classification. 
\end{proof}

For the computation of the second Chern class 
$c_2(N_{\widetilde S/\widetilde Y})$, we prepare an easy lemma.

\begin{lemma}\label{easy}
$\widetilde E^4=c_2(N_{\widetilde S/\widetilde Y})$.
\end{lemma}

\begin{proof}
Note that $\tilde f$ is a blowing-up of the smooth variety $\widetilde Y$ 
along the smooth centre $\widetilde S$. Therefore, we can regard 
$\tilde h=\tilde f _{|\widetilde E}$ as
$\tilde h: \widetilde E=\mathbb P(N_{\widetilde S/\widetilde X}) \to S$.  
We have a line bundle $L$ on $\widetilde E$ such that 
there is an exact sequence of bundle maps
\[
0\lto L\lto \tilde h^* N_{\widetilde S/\widetilde X}
\lto \mathscr{O}_{\widetilde E}(1)\lto 0 .
\]
By naturality, we have 
$c_2(\tilde h^* N_{\widetilde S/\widetilde X})
=c_2(N_{\widetilde S/\widetilde X})\cdot [F] \in H^4(\widetilde E,\mathbb{Z})$, 
where $F$ is a fibre of $\tilde h$. This implies
$c_2(N_{\widetilde S/\widetilde X})
=c_2(\tilde h^*N_{\widetilde S/\widetilde X})
\cdot \mathscr{O}_{\widetilde E}(1)
=L\cdot \mathscr{O}_{\widetilde E}(1)^2$. 
On the other hand, one has trivially 
$\tilde h^* c_1(N_{\widetilde S/\widetilde X})
=c_1(\tilde h^* N_{\widetilde S/\widetilde X}) 
= c_1(\mathscr O_{\widetilde E}(1))+c_1(L)$. 
Combining these we get
$c_2(N_{\widetilde S/\widetilde X})
= \tilde h^* c_1(N_{\widetilde S/\widetilde X})\cdot 
	c_1(\mathscr O_{\widetilde E}(1))^2 + \widetilde E^4$, 
for $\mathscr O_{\widetilde E}(1)=(-\widetilde E)_{|\widetilde E}$. 
But we have $c_1(N_{\widetilde S/\widetilde X})
=c_1(T_{\widetilde Y|\widetilde S}) - c_1(T_{\widetilde S})=0$ and 
therefore $c_2(N_{\widetilde S/\widetilde X})=\widetilde E^4$.
\end{proof}

\begin{remark}\label{conv}
We show the converse of Theorem \ref{main} using this lemma. 
Let $X=\Hilb^2 (S)$ for some K3 surface $S$.
Noting that \eqref{hilbdiagram} fits into the diagram \eqref{modeldiagram}, 
(i,ii) of the Theorem \ref{main} are evident. In the notation of 
\eqref{hilbdiagram} we have $N_{\Delta/S\times S}
\cong T_S\cong \Omega^1 _S$, for $S$ is K3. Therefore, 
for the exceptional divisor
$\widetilde E$ of $F$, we get $\widetilde E^4=c_2(N_{\Delta/S\times S})
=c_2(\Omega ^1_S)=24$. By the projection formula and ramification, we have 
\begin{equation}\label{EandEtild}
\widetilde E^4=\left(\frac 12 p^* E\right)^4
=\frac 1{16}(p^*E)^4=\frac 18E^4 ,
\end{equation}
so that $E^4=8\cdot 24=192$.
\end{remark}

\begin{corollary}\label{K3}
Under the assumption and notation as in {\rm \S 1}, $\widetilde S$ is
a K3 surface.
\end{corollary}

\begin{proof}
Assume the contrary, \emph{i.e.} assume $\widetilde S$ be
an abelian surface. Note that $\widetilde E^4=24$ by \eqref{EandEtild}. 
Lemmas \ref{appdhollef} and \ref{easy} imply $L_{hol}(g)
=\displaystyle \frac{24}{16}\not\in \mathbb Z$, which is a contradiction.

\end{proof}

Now is the time to prove Proposition \ref{quant}.

\begin{proof}[Proof of\; Proposition \ref{quant}]
\noindent (i) Let $\widehat E$, $\widehat{\widetilde E}$ be tubular
neighbourhoods of $E$, $\widetilde E$ and set
$X^{\circ}=X\backslash E$, $\widetilde X^{\circ}=\widetilde X
\backslash \widetilde E$. Then we have
\begin{eqnarray*}
X=X^{\circ}\cup \widehat E, \quad
	& X^{\circ}\cap \widehat E=\widehat E\backslash E,\\
\widetilde X=\widetilde X^{\circ}\cup \widehat{\widetilde E}, \quad
	&\widetilde X^{\circ}\cap \widehat{\widetilde E}
		=\widehat{\widetilde E}\backslash \widetilde E.
\end{eqnarray*}
Since the K3 surface $S\cong \widetilde S$ is simply connected
and $E\to S$ and $\widetilde E\to \widetilde S$ are 
$\mathbb P^1$-bundles, the homotopy exact sequence infers
\[
\pi _1(E)=\pi _1\left(\widehat E\right)=\{e\},\quad 
\pi _1\left(\widetilde E\right)
=\pi _1\left(\widehat{\widetilde E}\right)= \{e\}.
\]
Note that $\pi _1(X)=\{e\}$ because $X$ is irreducible symplectic.
By Van Kampen's theorem
\[
\xymatrix @ur {
\pi _1 (X^{\circ}) \ar[r] &\pi _1(X)\cong \{e\} \\
\pi _1\left(\widehat E\backslash E\right)\cong \mathbb Z
	\ar[u]_{\varphi} \ar[r]
	& \pi _1\left(\widehat E\right)\cong \{e\} \ar[u]
}\lower38pt\hbox{,}
\]
we know $\varphi$ is surjective. Consider the following commutative diagram
\[
\xymatrix{
\pi _1\left(\widetilde X^{\circ}\right)\ar[d]^{p_*}
	& \ar[l]^{\widetilde \varphi \quad}
	\pi _1\left(\widehat{\widetilde E}\backslash \widetilde E\right)
	\cong \mathbb Z \ar[d]_{p_*}\\
\pi _1(X^{\circ})
	&\ar[l]^{\varphi\quad}
	\pi _1\left(\widehat E\backslash E\right)\cong \mathbb Z
}\lower50pt\hbox{.}
\]
Since $\widetilde X^{\circ}\to X^{\circ}$ is \'etale of degree 2, 
$\pi _1(\widetilde X^{\circ})\overset{p_*}{\to}
\pi _1(X^{\circ})$ is injective but not surjective. It follows that 
$\widetilde \varphi$ is also surjective, because
$\pi _1\left(\widehat{\widetilde E}\backslash \widetilde E\right)
\to \pi _1\left(\widehat E\backslash E\right)$
is of index 2. Again using Van Kampen's theorem, we get
$\pi _1\left(\widetilde X\right)=\{e\}$, therefore $\pi _1
\left(\widetilde Y\right)=\{e\}$, for a birational map
$\tilde f$ does not change the fundamental group.

\noindent (ii) The condition $h^0(X,\Omega ^2_X)=1$ implies 
$h^0(Y,\Omega ^2_Y)=1$. Therefore, 
\[
L_{hol}(g)=1+\{1-(h^0(\widetilde Y,\Omega ^2_{\widetilde Y})-1)\}
	+1=4-h^0(\widetilde Y,\Omega ^2_{\widetilde Y}).
\]
On the other hand, by Lemmas \ref{appdhollef}, \ref{easy}, and
Corollary \ref{K3}, we see
\[
L_{hol}(g)=\frac 12+\frac{24}{16}=2 .
\]
Combining these, we get $h^0(\widetilde Y,\Omega ^2_{\widetilde Y})
=2$. 
\end{proof}

\section{Conclusion of the proof of Theorem \ref{main}}
\label{sec:3}

In this section, we complete the proof of Theorem \ref{main}. 

In the last section, we have shown that $\widetilde Y
\cong T_1\times T_2$, where, $T_1$, $T_2$ are K3 surfaces. 
To prove Theorem \ref{main}, we investigate the action of $g$ 
on $\widetilde Y$.

\begin{proposition}\label{P2}
$T_1\cong T_2\cong \widetilde S$. We may assume $\widetilde S$
is the diagonal of $\widetilde Y\cong \widetilde S\times \widetilde S$.
\end{proposition}

\begin{proof}
Let $p_i : \widetilde Y\to T_i\; (i=1,2)$ be the projections. 
To prove the proposition, it is enough to show
\begin{equation}\label{cl1}
\dim p_1(\widetilde S)=\dim p_2(\widetilde S)=2 ,
\end{equation}
for these imply $p_{i\; |\widetilde S}=:\varphi _i : 
\widetilde S\to T_i\; (i=1,2)$ is generically finite, but 
$K_{\widetilde S}\sim 0$ and $K_{T_i}\sim 0$ imply 
$\varphi _i$ has neither exceptional divisor nor ramification 
divisor, \emph{i.e.} $\varphi _i$ is isomorphism. 
In the following, we show \eqref{cl1} via case by case 
consideration.

If $\dim p_1(\widetilde S)=0$, \emph{i.e.} $p_1(\widetilde S)
=\{t\}\subset T_1$, we get $\widetilde S=p_1^{-1}(t)$, therefore 
$N_{\widetilde S/\widetilde Y}
=\mathscr{O}_{\widetilde S}^{\oplus 2}$ so that 
$c_2(N_{\widetilde S/\widetilde Y})=0$. This contradicts to 
Lemma \ref{easy}. By the same argument for $p_2$ 
we have $\dim p_i (\widetilde S)\geqslant 1\; (i=1,2)$.

\begin{claim}
If $\dim p_1(\widetilde S)=1$, then $\dim p_2 (\widetilde S)=2$. 
Moreover for $Z=p_1^{-1}(p_1(\widetilde S))$, we have $g(Z)=Z$.
\end{claim}

\begin{proof}[Proof of the claim]
Assume $\dim p_1(\widetilde S)=\dim p_2(\widetilde S)=1$. 
Let $C_i=p_i(\widetilde S)\; (i=1,2)$. Then we have
\[
\widetilde S\subset p_1^{-1} (C_1)\cap p_2^{-1}(C_2)
=C_1\times C_2.
\]
Since $\widetilde S$ is irreducible, so are $C_1$ and $C_2$.
But this implies that $\widetilde S$ is actually the product of 
these two curves, which is absurd.
Therefore we may assume $\dim p_2(\widetilde S)=2$ and $\varphi _2 : 
\widetilde S\overset{\sim}{\to} T_2$ is an isomorphism.
The second assertion of the claim is equivalent to $p_1(g(Z))=C_1$. 
Assume the contrary, \emph{i.e.} $p_1(g(Z))=T_1$. 
Then, we have the following diagram
\[
\xymatrix{
Z=C_1\times T_2 \ar[r]^{\quad g} \ar@{}[d]|{\bigcup} 
	\ar `l[d] `[dd]_{p_{1\; |Z}} [dd]
	& g(Z)\ar@{}[d]|{\bigcup} \ar@/^1pc/[ddr]^{\psi _1=p_{1\; |g(Z)}}\\
\widetilde S \ar[d]^{\varphi _1} \ar@{=}[r]
	& \widetilde S \ar[d]_{\varphi _1}\\
C_1 \ar@{=}[r] & C_1 \ar[r]^{\subset} & T_1
}\lower79pt\hbox{.}
\]
Clearly $\kappa (Z)=-\infty$, where $\kappa(Z)$ the Kodaira 
dimension of $Z$. The sub-additivity property of Kodaira dimension 
(cf. \cite{K}) implies that irreducible components of any fibre of 
$\psi _1$ are rational curves. 
On the other hand, connected components of the general fibre of
$\varphi _1 : \widetilde S\to C$ are elliptic curves so that  
$\psi _1$ contains an elliptic curve as its fibre, a contradiction.
\quad (End of the proof of the claim)
\end{proof}

To prove the proposition, what we have to do is to get a contradiction 
assuming $\dim p_1(\widetilde S)=1$ and $g(Z)=Z$, thanks to the claim.
Consider the diagram
\[
\xymatrix@R=10pt{
& Z=C_1\times T_2 \ar[dl]_{p_1}\ar@{}[d]|{\bigcup}\ar[dr]^{p_2}\\
C_1 & \widetilde S\ar[l]_{\varphi _1} \ar[r]^{\varphi _2}
& T_2
}\lower26pt\hbox{.}
\]
Let $U_t= p_1^{-1}(t)\; (t\in C_1)$ and consider the normalizations
\[
\xymatrix{
\overline{Z}=\mathbb P^1\times T_2\ar[d]^{\bar p_1} \ar[r]^{norm.}
& Z=C_1\times T_2 \ar[d]^{p_1}\\
\mathbb P^1\ar[r]^{norm.}& C_1
}\lower39pt\hbox{.}
\]
The involution $g$ on $Z$ ascends to $\overline Z$. Since $\bar p_1$ is 
the anti-canonical map of $\overline Z$, $g$ descends to $\mathbb P^1$. 
This implies $\dim p_1(g(U_t))=0$ therefore $p_1(g(U_t))=\{t\}$, since 
$U_t\cap \widetilde S=g(U_t)\cap \widetilde S\neq \emptyset$. 
Then we have a family 
of automorphisms $\{\psi _t : U_t\to U_t\}$. There is a commutative diagram 
of isomorphisms
\[
\xymatrix{
U_t \ar[d]_{\psi _t} \ar[r]^{p_{2\; |U_t}} & T_2 \ar[d]
	& \widetilde S \ar[l]_{p_{2\; |\widetilde S}} \ar[d]^{\rho _t}\\
U_t \ar[r]^{p_{2\; |U_t}} & T_2 
	& \widetilde S \ar[l]_{p_{2\; |\widetilde S}}
}\lower42pt\hbox{.}
\]
Since a K3 surface has no infinitesimal automorphism, $\rho _t$ is 
independent of $t$. But $\psi _t$ fixes the points on $D_t=U_t\cap 
\widetilde S$ and $\widetilde S=\bigcup _{t\in C_1} D_t$, $\rho _t$ 
induces the identity on $\widetilde S$ and also on $U_t$. This contradicts 
$g\neq \id$. 
\end{proof}

Finally the following proposition completes our proof of the Main 
Theorem, in view of Example \ref{hilb}.

\begin{proposition}
Notation as above. The action of $g$ on $\widetilde Y=\widetilde S\times
\widetilde S$ is the permutation of two components.
\end{proposition}

\begin{proof}
Consider the diagram
\[
\xymatrix{
\widetilde S \ar[r]^{\mbox{\footnotesize diag.}\qquad } 
	& \widetilde Y=\widetilde S \times \widetilde S 
	\ar[d]^{p_1} \ar[r]^{\qquad p_2} & \widetilde S\\
& \widetilde S &
}\lower44pt\hbox{.}
\]
Let $T_t=p_1^{-1} (t)\; (t\in \widetilde S)$ and 
$\Gamma _i=\{t\in \widetilde S\mid \dim p_i\circ g(T_t)=1\}$. 
Note that $\Gamma _i$ is a locally closed set in Zariski topology.

\begin{claim}
$\dim \Gamma _i\leqslant 1\quad (i=1,2)$.
\end{claim}

\begin{proof}[Proof of the claim]
Assume the contrary,  \emph{i.e.} assume that 
$\Gamma _1\subset \widetilde S$ contains an Zariski 
open set. For any $t\in \widetilde S$, 
\[
g(T_t)\overset{p_1}{\lto} C_t\subset \widetilde S
\] 
with $C_t$ is a rational curve containing $t$. Since $t$ sweeps over an 
open set  of $\widetilde S$, 
$\{C_t\}$ is a covering family of rational curves on $\widetilde S$. 
This is impossible because $\widetilde S$ is a K3 surface.
\quad (End of the proof of the claim)
\end{proof}

Therefore, there exists a Zariski open set $V\subset \widetilde S$ 
such that 
\[
(\dim p_1\circ g(T_t),\dim p_2\circ (T_t))
=(0,2)\mbox{, }(2,2)\mbox{ or }(2,0)
\]
for any $t\in V$. 

In the first case, we have $g(T_t)=T_t$ and $g$ induces 
$j_t\in \Aut (T_t)\cong \Aut (\widetilde S)$. 
Since $\widetilde S$ has no infinitesimal automorphism, $j_t$ is constant 
with respect to $t\in \widetilde S$. But $j_t(t)=t$ for any 
$t\in V$ implies $j_t=\id$, contradiction.

In the second case, $g(T_t)$ is the graph of an automorphism 
$f_t\in \Aut (\widetilde S)$. Since $g(T_t)\cap g(T_{t'})=\emptyset$ 
for $t\neq t'$, we have $f_t\neq f_{t'}$. This contradicts
the discreteness of $\Aut (\widetilde S)$. 

Therefore only the last case can happen. This implies 
$g(p_1^{-1}(t))=p_2^{-1}(t)$, for $T_t\cap \widetilde S\neq \emptyset$.
Exchanging the roles of $p_1$ and $p_2$, we also have 
$g(p_2^{-1}(t))=p_1^{-1}(t)$. Finally we get
\[
\xymatrix{
\{(s,t)\} \ar@{=}[r] & p_1^{-1}(s)\cap p_2^{-1}(t) \ar[d]^g\\
&p_2^{-1}(s)\cap p_1^{-1}(t) \ar@{=}[r] & \{(t,s)\}
}
\]
for $s,t\in V$. This shows $g$ is a permutation of $p_1$ and $p_2$.
\end{proof}

\paragraph{{\bf Acknowledgement}} 

The author would like to express his profound gratitude to 
Prof. Yujiro Kawamata, his supervisor, for comments and continuous 
encouragement. He also thanks Tetsushi Ito and Shunsuke Takagi for 
stimulating discussions, Prof. Keiji Oguiso for his comments 
and indicating an improvement of the proof of Proposition \ref{P2}. 
He is also indebted to the referee for his careful reading of 
the manuscript and many valuable comments and suggestions.

%%%%%%


\end{document}